\input amstex 
\documentstyle{amsppt}
\input bull-ppt
\keyedby{bull274/lic}
\define\spn{\operatorname{Span}}
\define\om{\overline\mu}

\topmatter
\cvol{26}
\cvolyear{1992}
\cmonth{April}
\cyear{1992}
\cvolno{2}
\cpgs{288-293}
\title The Mackey-Gleason Problem \endtitle
\author L. J. Bunce and J. D. Maitland Wright\endauthor
\shortauthor{L. J. Bunce and J. D. Maitland Wright}
\shorttitle{The Mackey-Gleason Problem}
\address Analysis and Combinatorics Research Centre, 
Mathematics Department, University of Reading,
P O Box 220, Whiteknights, Reading RG6 2AX, England 
\endaddress
\date May 21, 1991\enddate
\abstract Let $A$ be a von Neumann algebra with no direct 
summand of
Type $\roman I_2$, and let $\scr P(A)$ be its lattice of 
projections. 
Let $X$ be a Banach space. Let $m\:\scr P(A)\to X$ be a 
bounded function
such that $m(p+q)=m(p)+m(q)$ whenever $p$ and $q$ are 
orthogonal 
projections. The main theorem states that $m$ has a unique 
extension
to a bounded linear operator from $A$ to $X$. In 
particular, each bounded
complex-valued finitely additive quantum measure on $\scr 
P(A)$ has a unique
extension to a bounded linear functional on $A$.\endabstract
\endtopmatter

\document
\heading Physical background\endheading
In von Neumann's approach to the mathematical foundations 
of quantum
mechanics, the bounded observables of a physical system 
are identified
with a real linear space, $L$, of bounded selfadjoint 
operators on a
Hilbert space $H$. It is reasonable to assume that $L$ is 
closed in
the weak operator topology and that whenever $x\in L$ then 
$x^2\in L$.
(Thus $L$ is a Jordan algebra and contains spectral 
projections.) Then
the projections in $L$ form a complete orthomodular 
lattice, $\scr P$,
otherwise known as the lattice of ``questions'' or the 
quantum logic
of the physical system. A quantum measure is a map 
$\mu\:\scr P\to
\Bbb R$ such that whenever $p$ and $q$ are orthogonal 
projections $\mu
(p+q)=\mu(p)+\mu(q)$. 

In Mackey's formulation of quantum mechanics \cite{11} his 
Axiom VII
makes the assumption that $L=L(H)_{sa}$. Mackey states, 
that in 
contrast to his other axioms, Axiom VII has no physical 
justification; 
it is adopted for mathematical convenience. One of the 
technical advantages
of this axiom was that, by Gleason's Theorem, a completely 
additive
positive quantum measure on the projections of $L(H)$ is 
the restriction
of a bounded linear functional (provided $H$ is not 
two-dimensional).
In order to weaken Axiom VII it was desirable to 
strengthen Gleason's
Theorem.

\heading Introduction\endheading
Let $\scr P(A)$ be the lattice of projections in a von 
Neumann algebra $A$,
let $X$ be a Banach space, and let $\mu\:\scr P(A)\to X$ 
be a function 
such that 
\roster
\item "(a)" $\mu(e+f\,)=\mu(e)+\mu(f\,)$ whenever $ef=0$, 
\item "(b)" $\sup\{\|\mu(e)\|\: e\in\scr P(A)\}<\infty$. 
\endroster
Then $\mu$ is said to be a {\it finitely additive, 
$X$\<-valued measure
on $\scr P(A)$}.

Clearly each bounded linear operator from $A$ to $X$ 
restricts to a finitely
additive $X$\<-valued measure. When $A$ is the algebra of 
two-by-two
matrices and $X$ is one-dimensional, there exist examples 
of measures that
fail to extend to linear\ functionals.

Our main result is

\thm{Theorem A} Let $A$ be a von Neumann algebra with no 
direct summand
of Type $\roman I_2$. Then, for each Banach space $X$, 
each $X$\<-valued
measure on $\scr P(A)$ has a unique extension to a bounded 
linear operator
from $A$ to $X$.\ethm

This immediately specializes to give

\thm{Theorem B} Let $A$ be a von Neumann algebra with no 
direct summand
of Type $\roman I_2$. Then each complex-valued finitely 
additive measure
on $\scr P(A)$ extends to a bounded linear functional on 
$A$.
\ethm

In fact, we shall see that Theorem A follows easily from 
Theorem B.
We note, however, that to deduce Theorem A it is essential 
to have Theorem B
for all real-valued finitely additive measures. It does 
not suffice to
know this result for positive measures or for countably 
additive measures.
The lack of positivity causes considerable difficulty in 
establishing
this theorem.

Theorem B answers a natural question first posed by 
G.~W.~Mackey some
thirty years ago. When $\mu$ is positive, that is, when 
$\mu(e)\ge 0$
for each projection $e$. Theorem B was established by 
Christensen \cite7
for properly infinite algebras and algebras of Type 
$\roman I_n$ and by
Yeadon \cite{15, 16} for algebras of finite type. The 
first major progress
had been made by Gleason \cite9, who, by using an 
ingenious geometric
argument, settled the question for positive completely 
additive measures
on the projections of $L(H)$. (See \cite8 for an 
elementary proof of this
deep result.) Aarnes \cite1 and Gunson \cite{10} made 
important contributions,
especially concerning continuity properties. Paszkiewicz 
\cite{13}, working
independently of Christensen and Yeadon, established 
Theorem B for 
$\sigma$\<-finite factors, however, to extend his results 
to nonfactorial
von Neumann algebras, he requires $\mu$ to be positive and 
countably
additive. A lucid and meticulous exposition of Theorem B 
for positive
measures is given by Maeda \cite{12}.

\heading 1. Vector measures\endheading
The following short argument shows that Theorem A is a 
consequence of 
Theorem B.

\thm{Lemma 1.1} Let $A$ be a von Neumann algebra such that 
each finitely
additive \RM(complex\/\RM) measure on $\scr P(A)$ has an 
extension to a
bounded linear functional. Let $X$ be a Banach space and 
let $m\:\scr P
(A)\to X$ be a finitely additively $X$\<-valued measure. 
Then $m$ has a
unique extension to a bounded linear operator from $A$ to 
$X$.\ethm

\demo{Proof} Elementary estimates show that when $\beta\in 
A^*$ then
$$
\|\beta\|\le 4\sup\{|\beta(p)|\: p\in\scr P(A)\}\.
$$
Let $K$ be a constant such that $\|m(p)\| \le K$ for each 
$p\in\scr P(A)$.

For any $\phi \in X^*$, $p\to\phi m(p)$ is a 
complex-valued finitely
additive measure on $\scr P(A)$. By hypothesis there 
exists $\beta\in
A^*$ such that $\beta(p)=\phi m(p)$ for each $p\in\scr 
P(A)$. 

Let $x=\sum_1^n\lambda_j p_j$ be a finite linear 
combination of
projections $p_1\!,p_2\!,\dots, p_n$\<. Then
$$
\phi\lf(\sum_1^n\lambda_j m(p_j)\rt)=\sum_1^n\lambda_j\phi 
m(p_j)
=\beta\lf(\sum_1^n \lambda_j p_j\rt)=\beta(x)\.
$$
So
$$
\align
\lf|\phi\lf(\sum_1^n\lambda_j m(p_j)\rt)\rt|
\le & 4\|x\|\, \|\phi\|\sup\{\|m(p)\|\: p\in P(A)\}\\
\le & 4\|x\|\, \|\phi\| K\.
\endalign
$$
It follows from the Hahn-Banach Theorem that
$$
\lf\|\sum_1^n \lambda_j m(p_j)\rt\| \le 4K\|x\|\.
$$
In particular, $\sum_1^n \lambda_j p_j=0$ implies 
$\sum_1^n\lambda_j
m(p_j)=0$. So $m$ has a unique extension to a linear 
operator $T\:
\spn\scr P(A)\to X$. 
$$
\lf\| T\lf(\sum_1^n \lambda_j p_j\rt)\rt\|
= \lf\|\sum_1^n \lambda_j m(p_j)\rt\| \le 4K\lf\|\sum_1^n 
\lambda_j p_j\rt\|\.
$$
So $T$ is bounded and hence has a unique extension to a 
bounded linear
operator from $A$ to $X$. \enddemo

\heading 2. Scalar measures\endheading
In all that follows $A$ is a von Neumann algebra with no 
direct summand
of Type $\roman I_2$. Let $\scr P(A)$ be the lattice of 
projections in
$A$.

Since each complex-valued measure on $\scr P(A)$ is of the 
form $\mu+i\nu$,
where $\mu$ and $\nu$ are real-valued measures, it 
suffices to prove Theorem 
B for real-valued measures.

From now onward, $\mu$ is a finitely additive real-valued 
measure on
$\scr P(A)$. That is, $\mu\: \scr P(A)\to\Bbb R$ is a 
function such that
\roster
\item "(a)" $\mu(p+q)=\mu(p)+\mu(q)$ whenever $p$ and $q$ 
are orthogonal
projections;
\item "(b)" $\sup\{|\mu(p)|\: p\in\scr P(A)\}<+\infty$. 
\endroster
We define the variation of $\mu$, $V$, by 
$$
V(p)=\sup\{|\mu(e)|\: e\le p\};
$$
we also define
$$
\alpha(p)=\sup\{\mu(e)\: e\le p\}\.
$$

Straightforward arguments show that $\mu$ has a unique 
extension to a
function $\om\: A\to\Bbb C$, where $\om$ is linear and 
bounded on each
abelian $^*$\<-subalgebra of $A$ and where $\om(x+
iy)=\om(x)+i\om(y)$
whenever $x$ and $y$ are selfadjoint. Moreover, it can be 
shown \cite4 that
\roster
\item $\sup\{|\om(x)|\: x=x^*$ and 
$\|x\|\le1\}=2\alpha(1)-\mu(1)$,
\item $\sup\{\om(x)\: 0\le x\le 1\}=\alpha(1)$.
\endroster

In order to see (2), let $0\le x\le 1$. Then, for suitable 
spectral projections,
$x=\sum_1^\infty 2^{-n} e_n$. Hence, for some $n$, 
$\mu(e_n)\ge\om(x)$. 

The lack of positivity of $\mu$ greatly increases the 
difficulties of
establishing the linearity of $\om$ when $A$ is properly 
infinite or when
$A$ is of Type $\roman{II}_1$. However, when $A$ is of 
Type $\roman I_n$
$(n\ne 2)$ linearity can be established by a fairly 
straightforward extension
of the arguments for positive measures.

The first step is to notice that when $\mu$ is a measure 
on the projections
of $M_n(\Bbb C)$, the algebra of $n\times n$ matrices over 
$\Bbb C$, then
when $T$ is the canonical (unnormalized) trace on 
$M_n(\Bbb C)$ we have
$\mu(e)\le\alpha(1)=\alpha(1) T(e)$ for each minimal 
projection $e$. So
$\alpha(1) T-\mu$ is a positive measure on $\scr 
P(M_n(\Bbb C))$. Hence,
by Gleason's Theorem for finite-dimensional Hilbert 
spaces, $\alpha(1)
T-\om$ is linear (provided $n\ne 2$\<). We now revert to 
the general
situation and conclude that $\om$ is linear on each 
subalgebra of $A$
that can be embedded in a subalgebra of $A$ that is 
isomorphic to a 
Type $\roman I_n$ factor $(n\ge 3)$. By elementary 
algebraic arguments,
it can be shown \cite8 that if $B$ is a subalgebra of $A$ 
and $B\approx
M_2(\Bbb C)\subset\Bbb D$, then either $B\subset C\subset 
A$, where
$C\approx M_4(\Bbb C)$ or $B\subset C\oplus D$, where 
$C\simeq M_4
(\Bbb C)$ and $D\subset E\subset A$ with $E\approx 
M_3(\Bbb C)$. Hence
$\om$ is linear on $B$.

By ``patching'' together Type $\roman I_2$ factorial 
subalgebras of $A$
it can be shown \cite2, building on techniques of 
Christensen \cite7,
that $\mu$ is uniformly continuous on $\scr P(A)$ and 
$\om$ is linear
on each Type $\roman I_n$ subalgebra of $A$. In 
particular, $\om$ is
linear on $W(1,p,q)$, where $W(1,p,q)$ is the 
$W^*$\<-subalgebra of $A$
generated by the identity and an arbitrary pair of 
projections $p$ and
$q$. 

\heading 3. Approximate linearity\endheading
The restriction of $\om$ to the centre of $A$ is linear. 
So, by the
Hahn-Banach Theorem, there exists $\sigma \in A^*$ such 
that $\om$
and $\sigma$ coincide on the centre. By replacing $\om$ by 
$\om-\sigma$
if necessary, we can assume that $\om$ vanishes on the 
centre of $A$.
By dividing by a suitable constant we can also assume that 
$\sup\{|\overline\mu
(x)|\: x=\lambda^*, \|\lambda\|\le 1\}=1$.
It can then be shown that $\alpha(1)=\tfrac 12$. 

Lack of positivity leads to a number of difficulties 
establishing

\thm{Lemma 3.1} Let $A$ be properly infinite. Let 
$0<\delta<\tfrac 12$.
There exists a projection $e$ in $A$ with $1\sim e\sim 
1-e$ and 
such that $\tfrac 12-\delta^2<\mu(e)$. Then, for each 
$p\in\scr P(A)$. 
$$
|\mu p-\om(epe)-\om (1-e) p(1-e)| <5\delta\.
$$
\ethm

The next lemma depends on the fact that $\om$ is linear on 
each 
$W^*$\<-subalgebra generated by a pair of projections.

\thm{Lemma 3.2} Let $0<\varepsilon<1$ and let $m$ be such 
that
$\sum_{m+1}^\infty 2^{-n}<\varepsilon$. Let $e$ be a 
projection such
that

\RM{(1)} $|\mu(p)-\om(epe)-\om((1-e)p(1-e))|<\varepsilon/m$
for each $p\in\scr P(A)$\,\RM;

\RM{(2)} $|\om(a+b)-\om(a)-\om(b)|<\varepsilon/m$
whenever $a\ge 0$, $b\ge 0$, and $a+b\le e$\,\RM;

\RM{(3)} $|\om(c+d)-\om(c)-\om(d)|<\varepsilon/m$ 
whenever $c\ge 0$, $d\ge 0$, and $c+d\le 1-e$. 
Then, whenever $x\ge 0$, $y\ge 0$, and $x+y\le 1$, 
$|\om(x+y)-\om(x)-\om(y)|<20\varepsilon$. \ethm

We shall now sketch a proof of Theorem B for $A$ a 
properly infinite
von Neumann algebra.

Let $a$ and $b$ be fixed, positive elements of $A$ with $a+
b\le 1$.
Choose $\varepsilon$ with $0<\varepsilon<1$. Let $m$ be 
such that
$\sum_{m+1}^\infty 2^{-n}<\varepsilon$. We put 
$\delta=\varepsilon/
42m$. Let $e$ be a projection that satisfies the 
conditions of Lemma
3.1. 

Let $x\ge 0$, $y\ge 0$ with $x+y\le e$. 

We find a projection $f\le 1-e$ such that $f\sim 1$ and 
$V(f\,)<\delta^2$.
Put $h=e+f$. We then find three orthogonal projections, 
majorized by $f$ and
each equivalent to $e$. By applying the $4\times 4$ matrix 
construction
due to Christensen \cite7 we can find orthogonal 
projections $p$ and $q$,
majorized by $h=e+f$, with $x=2epe$ and $y=2eqe$. Since 
$x$ is in the
$^*$\<-algebra generated by $e$ and $p$,
$$
|2\mu(p)-\om(x)| =|\om(2p-x)|=2|\om(p-epe)|\.
$$
Since $e$, $f$, and $p$ are in the properly infinite 
$W^*$\<-algebra
$hAh$ and $e+f=h$, we see that $e$ is in the 
$W^*$\<-subalgebra of $hAh$
generated by $f$ and $p$. So
$$
2|\om(p-epe)|\le 2|\om (fpf\,)|+2 |\om(fpe+epf\,)|\.
$$
By considering $\om$ restricted to $fAf$ we can show that
$$
|\mu(fpf\,)|\le 2V(f\,)<2\delta^2\.
$$
Also, 
$$
fpe+epf=(1-e) pe+ep(1-e)\.
$$
So, applying Lemma 3.1, $|\om(fpe+epf\,)|<5\delta$. So
$$
|2\mu(p)-\om(x)| <2\delta^2+10\delta<14\delta\.
$$
Similarly,
$$
|2\mu(q)-\om(y)|< 14\delta\quad\text{and}\quad|2\mu(p+
q)-\om(x+y)|<14\delta\.
$$
So
$$
|\om(x+y)-\om (x)-\om(y)|<42\delta=\tfrac\varepsilon m\.
$$

We may repeat the above argument, interchanging the roles 
of $e$ and
$1-e$, to deduce that whenever $z\ge 0$, $w\ge 0$, and $z+
w<1-e$,
$$
|\om(w+z)-\om(w)-\om(z)|<\tfrac \varepsilon m\.
$$
We now appeal to Lemma 3.2 to obtain
$$
|\om(a+b)-\om(a)-\om(b)|<20\varepsilon\.
$$
Since $\varepsilon$ is arbitrary, it follows that
$$
\om(a+b)=\om(a)+\om(b)\.
$$
Hence $\om$ is linear.

When $\mu$ is $\sigma$\<-additive, we can give a 
reasonably straightforward
proof of linearity for Type $\roman{II}_1$ algebras. 
However, in order to
obtain Theorem A, it is essential to obtain Theorem B when 
$\mu$ is finitely
additive, not positive and not $\sigma$\<-additive. For 
von Neumann algebras
of finite type this forces us to use a more elaborate 
argument.

\heading 4. Open problem\endheading
Let $A$ be a (unital) $C^*$\<-algebra. Let $\phi\: 
A\to\Bbb C$ be a
function whose restriction to each abelian 
$^*$\<-subalgebra is linear,
is such that $\{|\phi(x)|\: \|x\|\le 1\}$ is bounded and 
whenever $x$
and $y$ are selfadjoint, $\phi(x+iy)=\phi(x)+i\phi(y)$. 
Then $\phi$ is
said to be a quasi-linear functional. 

\ex{Problem} For which $C^*$\<-algebras $A$ is it true 
that every
quasi-linear functional on $A$ is linear?\endex

We have the following consequence of Theorem B. 

\thm{Corollary} Let $M$ be a von Neumann algebra with no 
direct summand
of Type $\roman I_2$. Let $I$ be a closed ideal of $M$ and 
let $A=M/I$.
Then every quasi-linear functional on $A$ is linear.\ethm

\demo{Proof} Let $\phi\: A\to\Bbb C$ be quasi-linear. Let 
$\pi\: M\to M/I$
be the canonical quotient homomorphism. Then the 
restriction of $\phi\pi$
to $\scr P(M)$ is a finitely additive measure that, by 
Theorem B, has a
unique extension to a bounded linear functional on $M$. 
Hence $\phi$ is
linear. \enddemo

\thm{Corollary} Each quasi-linear functional on the Calkin 
algebra is
linear.\ethm

\enddocument